\documentclass[12pt]{amsart}
\usepackage{amssymb, mathrsfs,amsmath}
\usepackage{latexsym}
\usepackage{amsfonts}
\usepackage{shadow}
\date{}
\author[D. Alpay]{Daniel Alpay}
\address{(DA) Department of Mathematics \newline
Ben Gurion University of the Negev \newline P.O.B. 653, \newline
Be'er Sheva 84105, \newline ISRAEL}
\email{dany@math.bgu.ac.il}
\author[H. Attia]{Haim Attia}
\address{(HA) Department of Mathematics \newline
Sami Shamoon College of Engineering
 \newline
Be'er Sheva 84100, \newline ISRAEL} \email{haima@sce.ac.il}
\keywords{White noise space, linear systems on commutative rings, interpolation problems}
\subjclass{Primary: 60H40, 93E10. Secondary: 30A80}

\thanks{D. Alpay thanks the
Earl Katz family for endowing the chair which supported his
research}
\title[Interpolation problem  in a
commutative ring]{An interpolation problem for functions with values in a
commutative ring}
\dedicatory{Dedicated to the memory of Israel Gohberg}
\begin{document}
\keywords{white noise space, stochastic distributions, linear systems on rings}
\subjclass{60H40, 93C05}
\maketitle
\parindent 0cm
\newtheorem{Pa}{Paper}[section]
\newtheorem{Tm}[Pa]{{\bf Theorem}}
\newtheorem{La}[Pa]{{\bf Lemma}}
\newtheorem{Cy}[Pa]{{\bf Corollary}}
\newtheorem{Rk}[Pa]{{\bf Remark}}
\newtheorem{Pn}[Pa]{{\bf Proposition}}
\newtheorem{Pb}[Pa]{{\bf Problem}}
\newtheorem{Dn}[Pa]{{\bf Definition}}
\newtheorem{Ex}[Pa]{{\bf Example}}
\renewcommand{\theequation}{\thesection.\arabic{equation}}
\def\L{\mathbf L}
\def\R{\mathbb R}
\def\N{\mathbb N}
\def\C{\mathbb C}
\def\D{\mathbb D}
\def\s{\mathscr S(\R)}
\def\ss{\mathscr S'(\R)}
\def\(s){\mathscr S(\R\times\R)}
\def\F{\mathcal F}
\def\P{\mathcal P}
\def\W{\mathcal W}
\def\Dom{{\rm dom}~(T_m)}
\def\Doms{{\rm dom}~(T_m^*)}
\def\Def{\stackrel{{\rm def.}}{=}}
\def\rr{\mathbf r}
\def\w{\omega}
\def\E{\mathscr I}
\begin{abstract}
It was recently shown that the theory of linear stochastic systems can be
viewed as a particular case of the theory of linear systems on a certain
commutative ring of power series in a countable number of variables. 
In
the present work we study an interpolation problem in  this setting. A key tool
is the principle of permanence of algebraic identities.

\end{abstract}
\keywords{}
\subjclass{}
\maketitle
\tableofcontents
\section{Introduction}
There are numerous connections between classical interpolation problems
and optimal control and the theory of linear systems; see for instance
\cite{bgr, MR2002b:47144}. In these settings, the coefficient space
is the complex field $\mathbb C$, or in the case of real systems,
the real numbers $\mathbb R$. 
Furthermore, already from its inception,  linear system theory was considered
when the coefficient space is a general  (commutative) field, or more generally
a commutative ring; see \cite{MR0255260,SSR}.
In \cite{MR2610579,alp} a new approach to
the theory of
 linear stochastic systems was developed, in which the coefficient space is now a
 certain commutative ring $\mathfrak R$ (see Section \ref{sec:r} below).
The results from \cite{MR0255260,SSR} do not seem to be directly applicable to
this theory, and the specific properties of $\mathfrak R$ played a key role
in the argumemts in \cite{MR2610579,alp}.\\

The purpose of this work is to discuss the counterparts of
classical interpolation problems in this new setting.
To set the problems into perspective, we
begin this introduction with a short discussion of the
deterministic case. In the classical theory of linear systems,
input-output relations of the form
\begin{equation}
\label{eqinput}
 y_n=\sum_{k=0}^n h_{n-k}u_k,\quad n=0,1,\ldots
\end{equation}
where $(u_n)_{n\in\mathbb N_0}$ is called the input sequence,
$(y_n)_{n\in\mathbb N_0}$ is the output sequence, and
$(h_n)_{n\in\mathbb N_0}$ is the impulse response, play an
important role. The sequence $(h_n)_{n\in\mathbb N_0}$ may consist
of matrices (of common dimensions), and then the
input and output sequences consist of vectors of appropriate
dimensions. Similarly state space equations
\[
\begin{split}
x_{n+1}&=Ax_n+Bu_n,\\
y_n&=Cx_n+Du_n,\quad n=0,1,\ldots
\end{split}
\]
play an important role. Here $x_n$ denotes the state at time $n$,
and $A,B,C$ and $D$ are matrices with complex entries. The
transfer function of the system is
\[
h(\lambda)=\sum_{n=0}^\infty h_n{\lambda}^n,
\]
in the case \eqref{eqinput}, and
\[
h({\lambda})=D+{\lambda}C(I-{\lambda}A)^{-1}B
\]
in the case of state space equations, when assuming the state at
$n=0$ to be equal to $0$. Classical interpolation problems bear
various applications to the corresponding linear systems. See for
instance \cite[Part VI]{bgr}, \cite{kaash-field}. To fix ideas, we consider the case of
bitangential interpolation problem for matrix-valued functions
analytic and contractive in the open unit disk (Schur functions),
and will even consider only the Nevanlinna-Pick interpolation
problem in the sequel to keep notation simple, but it will be
clear that the discussion extends to more general cases. Recall (see \cite[\S 18.5 p. 409]{bgr})
that the bitangential interpolation problem may be
defined in terms of a septuple of matrices
$\omega=(C_+,C_-,A_\pi,A_\zeta,B_+,B_-,\Gamma)$ by the conditions
\[
\begin{split}
\sum_{{\lambda}_0\in\mathbb D}{\rm
Res}_{{\lambda}={\lambda}_0}({\lambda}I-A_\zeta)^{-1}B_+S({\lambda})&=-B_-,\\
\sum_{{\lambda}_0\in\mathbb D}{\rm
Res}_{{\lambda}={\lambda}_0}S({\lambda})C_-({\lambda}I-A_\pi)^{-1}&=C_+,\\
\sum_{{\lambda}_0\in\mathbb D}{\rm
Res}_{{\lambda}={\lambda}_0}({\lambda}I-A_\zeta)^{-1}B_+S({\lambda})C_-({\lambda}I-A_\pi)^{-1}&=\Gamma,
\end{split}
\]
where $A_\zeta$ and $A_\pi$ have their spectra in the open unit disk,
where $(A_\zeta,B_+)$ is a full range pair (that is, controllable) and
where $(C_-,A_\pi)$ is a null kernel pair (that is, observable).
We send the reader to \cite{bgr} for the definitions.
Moreover, $\Gamma$ satisfies moreover the compatibility condition
\[
\Gamma A_\pi-A_\zeta \Gamma=B_+C_++B_-C_-.
\]
Let $P$ be the matrix (see \cite[p. 458]{bgr})
\begin{equation}
P=\begin{pmatrix}P_1&\Gamma^*\\
\Gamma&P_2\end{pmatrix},
\end{equation}
where $P_1$ and $P_2$ are the solutions of the Stein equations
\[
\begin{split}
P_1-A_\pi^*P_1A_\pi&=C_-^*C_--C_+^*C_+\\
P_2-A_\zeta P_2A_\zeta^*&=B_+B_+^*-B_-B_-^*.
\end{split}
\]
Furthermore, and assuming the unknown function $S$ to be $\mathbb C^{p\times q}$-valued,
\[
J=\begin{pmatrix}I_p&0\\0&-I_q\end{pmatrix}.
\]
When $P$ is strictly positive, the
solutions of the interpolation problem are given in terms of a
linear fractional transformation based on a $J$-inner rational
function $\Theta$ built from the septuple $\omega$ via the formula (see \cite[(18.5.6) p. 410]{bgr})
\begin{equation}
\label{eq:theta}
\begin{split}
\hspace{-2cm}\Theta(\lambda)=
I+(\lambda-\lambda_0)\begin{pmatrix}C_+&-B_+^*\\ C_-&B_-^*\end{pmatrix}\begin{pmatrix}(\lambda I-A_\pi)^{-1}&0\\
0&(I-\lambda A_\zeta^*)^{-1}\end{pmatrix}\times\\
&\hspace{-10cm}\times P^{-1}\begin{pmatrix}
(I-\lambda_0A_\pi^*)^{-1}C_+^*&-(I-\lambda_0A_\pi^*)^{-1}C_-^*\\
(A_\zeta-{\lambda}_0I)^{-1}B_+&(A_\zeta-{\lambda}_0 I)^{-1}B_-\end{pmatrix}
\end{split}
\end{equation}
where ${\lambda}_0$ is fixed on the unit circle and such that the various inverses exist in the above formula.
An important
fact is that the entries of $P_1$ and $P_2$ are rational functions of the entries of the matrices of $\omega$.
As a consequence,
there exists a rational function $f(\lambda)$, built
from $\omega$ and such that the entries of $\Theta$ are
polynomials in $\lambda$, with coefficients which are themselves
polynomials in the entries of the matrices of $\omega$ with coefficients in $\mathbb Z$. This fact
will allow us in the sequel to use the principle of permanence of
identities (see \cite[p. 456]{MR1129886}), to extend
interpolation problem to a more general setting.\\

Allowing in \eqref{eqinput} the input sequence
$(u_n)_{n\in\mathbb N_0}$
 to consist of random variables has been
considered for a long time. On the
other hand, allowing also the
impulse response of the system to carry some randomness seems much more
difficult to tackle. Recently a new approach to the theory of
linear stochastic systems was developed using Hida's white noise
space theory \cite{MR1244577}, \cite{MR2444857},
\cite{MR1387829}, and Kondratiev's spaces of stochastic test
functions and distributions \cite{MR1408433}. In this approach,
see \cite{aal2}, \cite{al_acap}, \cite{alp}, the complex numbers
are replaced by random variables in the white noise space, or more
generally, by stochastic distributions in the Kondratiev space,
and the product of complex numbers is replaced by the Wick
product. For instance, \eqref{eqinput} now becomes
\begin{equation}
\label{eqinput1}
 y_n=\sum_{k=0}^n h_{n-k}\lozenge u_k,\quad n=0,1,\ldots
\end{equation}
where the various quantities are in the white noise space, or
more generally in the Kondratiev's space of stochastic
distributions, and $\lozenge$ denotes the Wick product. An
important role in this theory is played by a ring $\mathfrak R$
of power series in countably many variables with coefficients in
$\mathbb C$. This ring is endowed with a topology, which is that
of the dual of a countably normed nuclear space. See Sections
\ref{wns} and \ref{sec:r}. Let us denote by
\begin{equation}
\label{fps} \rr(z)=\sum_{\alpha\in\ell} r_\alpha z^\alpha,
\end{equation}
an element of $\mathfrak R$, where $\ell$ denotes the set of
sequences $(\alpha_1,\alpha_2,\ldots)$, whose entries are in
\begin{equation*}
\mathbb {N}_0=\left\{ 0,1,2,3,\ldots\right\},
\end{equation*}
and for which $\alpha_k\not =0$ for only a finite number of
indices $k$, and where we have used the multi-index notation
\[
z^\alpha=z_1^{\alpha_1}z_2^{\alpha_2}\cdots\quad
%\alpha!=\alpha_1!\alpha_2!\cdots,\quad{\rm with}\quad
\alpha\in\ell.
\]
The ring $\mathfrak R$ has the following properties:\\
$(P_1)$ If $\rr\in\mathfrak R\,\,{\rm and}\,\,
\rr (0,0,0,\cdots)\not=0$, then $\rr$ has an inverse in $\mathfrak R$.\\
$(P_2)$ If $\rr\in\mathfrak R^{n\times n}$ is such that
$\rr(0,0,0,\cdots) =0_{n\times n}$ and if $f$ is a function of
one complex variable, analytic in a neighborhood of the origin,
with Taylor expansion
\[
f(\lambda)=\sum_{p=0}^\infty f_p\lambda^p,
\]
then, the series
\[
f(\rr)\stackrel{\rm def.}{=}\sum_{p=0}^\infty f_p\rr^p
\]
converges in $\mathfrak R^{n\times n}$. Furthermore, if $g$ is
another function of one complex variable, analytic in a
neighborhood of the origin, we have
\begin{equation}
\label{multi} (fg)(\rr)=f(\rr)g(\rr).
\end{equation}
$(P_3)$ If $\rr(z)=\sum_{\alpha\in\ell} r_\alpha
z^\alpha\in\mathfrak R$, then $\rr^*(z)\stackrel{\rm
def.}{=}\sum_{\alpha\in\ell} r_\alpha^*
z^\alpha\in\mathfrak R$, where $r_\alpha^*$ denotes the conjugate of the complex number $r_\alpha$.\\

Property $(P_1)$ implies in particular that a matrix $\mathbf
{A}\in\mathfrak R^{n\times n}$ is invertible in $\mathfrak
R^{n\times n}$ if and only if $\det \mathbf A(0)\not=0$. This
fact, together with $(P_2)$, allows to define expressions such as
\begin{equation}
\label{real}
\mathscr H(\lambda)=\mathbf{D}+\lambda \mathbf { C}(I_n-
\lambda\mathbf{A})^{-1}\mathbf{B}=\mathbf
D+\sum_{k=1}^\infty\lambda^k \mathbf C\mathbf A^{k-1}\mathbf B,
\end{equation}
where $\mathbf A$, $\mathbf B$, $\mathbf C$, and $\mathbf D$ are
matrices of appropriate dimensions and with entries in $\mathfrak
R$, and where $\lambda$ is an independent complex variable. As
explained in \cite{alp} this is the transfer function of some
underlying linear systems, and is a rational function with
coefficients
in $\mathfrak R$.\\

The purpose of this paper is to explain how to tackle in the
present setting counterparts of some classical interpolation
problems which appear in the theory of linear systems. To
illustrate our strategy, we focus on the Nevanlinna-Pick
interpolation problem, but our method works the same for the
general bitangential interpolation problem.
%
%We use the approach to interpolation based on the de Branges
%Rovnyak spaces. See \cite{MR90g:47003}, \cite{dym-ot41},
%\cite{abds2}, \cite{abds3}, and see \cite{dbr1}, \cite{dbr2} for
%information on de Branges Rovnyak spaces.
The computations done in the classical theory (that is, when the
coefficient space consists of the complex numbers) extend to the
case where $\mathbb C$ is replaced by the ring $\mathfrak R$.
In some cases, such as Nevanlinna-Pick interpolation, this can be shown by direct computations. In the general case, one
needs to use the principle of permanence of identities, see \cite[p.
456]{MR1129886}. We note that there
are other commutative rings with properties $(P_1),(P_2)$ and
$(P_3)$ for which the above analysis is applicable. See
\cite{a_salomon}.\\

The paper consists of six sections besides the present
introduction. In the second section we review Hida's white noise
space setting and the Kondratiev spaces of stochastic
distributions. The definition and main properties of the ring
$\mathfrak R$ are given in Section \ref{sec:r}. In Section
\ref{sec:3} we define and study analytic functions from an open
set of $\mathbb C$ with values in $\mathfrak R$. In Section
\ref{sec:np} we consider the Nevanlinna-Pick interpolation
problem. In the last section we discuss the bitangential interpolation problem.
% and the Carath\'eodory-Fej\'er problem respectively.
%We conclude the paper with various remarks gathered in Section
%\ref{sec:rem}.
\section{The white noise space}
\label{wns} We here review Hida's white noise space theory and
the associated spaces of stochastic distributions introduced by
Kondratiev. See \cite{MR1244577}, \cite{MR2444857},
\cite{MR1408433}, \cite{MR1387829}. Let $\mathscr S(\mathbb R)$
denote the Schwartz space of smooth {\sl real-valued} rapidly
decreasing functions. It is a nuclear space, and by the
Bochner-Minlos theorem (see \cite[Th\'eor\`eme 2, p.
342]{MR35:7123}), there exists a probability measure on the Borel
sets $\mathcal B$ of the dual space $\mathscr S(\mathbb
R)^\prime\stackrel{\rm def.}{=}\Omega$ such that
\begin{equation}
\label{voltaire} e^{-\frac{\|s\|_{{\mathbf L}_2(\mathbb
R)}^2}{2}}=\int_{\Omega}e^{i \langle \w,s\rangle}dP(\w),\quad
\forall s\in\mathscr S(\mathbb R),
\end{equation}
where the brackets $\langle \cdot,\cdot\rangle$ denote the
duality between $\s$ and $\ss$. The probability space
\[
\W= (\Omega,\mathcal B,dP)
\]
is called the {\sl white noise probability space}. We will be
interested in particular in $\mathbf L_2(\W)$, called the {\sl
white noise space}. For $s\in\mathscr S(\mathbb R)$, let $Q_s$
denote the random variable
\[
Q_s(\w)=\langle \w, s\rangle.
\]
It follows from \eqref{voltaire} that
\[
\|s\|_{{\mathbf L}_2(\mathbb R)}=\|Q_s\|_{\mathbf L_2(\W)}.
\]
Therefore, $Q_s$ extends continuously to an isometry from
${\mathbf L}_2(\mathbb R)$ into $\mathbf L_2(\W)$. In the
presentation of the Gelfand triple associated to the white noise
space which we will use, we follow \cite{MR1408433}. The white
noise space $\mathbf L_2(\W)$ admits a special orthogonal basis
$(H_{\alpha})_{\alpha\in\ell}$,  indexed by the set $\ell$ and
built in terms of the Hermite functions $\widetilde{h}_k$ and of
the Hermite polynomials $h_k$ defined by
\[
H_\alpha(\w)=\prod_{k=1}^\infty
h_{\alpha_k}(Q_{\widetilde{h}_k}(\w)).
\]
 We  refer to \cite[Definition 2.2.1 p.
19]{MR1408433} for more information. In terms of this basis, any
element of $\mathbf L_2(\W)$ can be written as
\begin{equation}
\label{fps1} F=\sum_{\alpha\in\ell}f_\alpha H_\alpha,\quad
f_\alpha\in\C,
\end{equation}
with
\[
\|F\|_{\W}^2=\sum_{\alpha\in\ell}|f_\alpha|^2\alpha!<\infty.
\]
There are quite a number of Gelfand triples associated to
$\mathbf L_2(\W)$. In our previous works \cite{MR2414165},
\cite{al_acap}, and in the present one, we focus on the one
consisting of the Kondratiev space $S_1$ of stochastic test
functions, of $\W$, and of the Kondratiev space $S_{-1}$ of
stochastic distributions. To define these spaces we first
introduce for $k\in{\mathbb N}$ the Hilbert space ${\mathcal
H}_{k}$ which consists of series of the form \eqref{fps} such that
\begin{equation}
\label{michelle}
\|F\|_{k}\stackrel{\rm def.}{=}
 \left(\sum_{\alpha\in\ell}(\alpha!)^2|f_\alpha|^2
(2{\mathbb N})^{k\alpha}\right)^{1/2}<\infty,
\end{equation}
and the Hilbert spaces $\mathcal H^\prime_k$ consisting of
sequences $(f_\alpha)_{\alpha\in\ell}$ such that
\[
\|F\|^\prime_{k}\stackrel{\rm def.}{=}
 \left(\sum_{\alpha\in\ell} |f_\alpha|^2
(2{\mathbb N})^{-k\alpha}\right)^{1/2}<\infty.
\]
We note that, for $F\in\mathcal H_{k}^\prime$ we have
\begin{equation}
\label{lim=0} \lim_{\substack{p\ge k \\ p\rightarrow
\infty}}\|F\|^\prime_{p}=|f_{(0,0,0\ldots)}|^2,
\end{equation}
as can be seen, for instance, by applying the dominated
convergence theorem to an appropriate discrete measure. Following
the usage in the literature, we will also write the elements of
$\mathcal H^\prime_k$ as formal power series
$\sum_{\alpha\in\ell}f_\alpha H_\alpha$.
Note that $(\mathcal H_k)_{k\in\mathbb N}$ forms a decreasing sequence
of Hilbert spaces, with increasing norms, while $(\mathcal H_k^\prime)_{k\in\mathbb N}$ forms an increasing sequence
of Hilbert spaces, with decreasing norms. The spaces $S_1$ and
$S_{-1}$ are defined by the corresponding projective and inductive limits
\[
S_1=\bigcap_{k=1}^\infty \mathcal H_k\quad{\rm and}\quad
S_{-1}=\bigcup_{k=1}^\infty \mathcal H^\prime_{k}.
\]
%The spaces $S_1$ and $S_{-1}$ are nuclear.\\
The Wick product is defined on the basis
$(H_\alpha)_{\alpha\in\ell}$ by
\[
H_\alpha\lozenge H_\beta=H_{\alpha+\beta}.
\]
It extends to an everywhere defined and continuous map 
from $S_1\times S_1$ into itself
and from $S_{-1}\times S_{-1}$ into itself\footnote{The continuity
properties are proved in \cite{a_salomon} for a more general family of rings.}
Let $l>0$,  and let
$k>l+1$. Consider $h\in {\mathcal H}^\prime_{l}$ and $u\in
{\mathcal H}^\prime_{k}$. Then, V\r{a}ge's inequality holds:
\begin{equation}
\label{vage}
\|h\lozenge u\|^\prime_{k}\le
A(k-l)\|h\|^\prime_{l}\|u\|^\prime_{k},
\end{equation}
where
\begin{equation}
\label{vage111} A(k-l)=\left(\sum_{\alpha\in\ell}(2{\mathbb
N})^{(l-k)\alpha}\right)^{1/2}<\infty.
\end{equation}
See \cite[Proposition 3.3.2 p. 118]{MR1408433}. The following
result is a direct consequence of \eqref{vage111} and will be
useful in the sequel.
\begin{La}
Let $F\in\mathcal H_p^\prime$. Then, $F^{\lozenge n}\stackrel{\rm
def.}{=}\underbrace{F\lozenge\cdots\lozenge F}_{n\,\,{times}}\in
\mathcal H_{p+2}^\prime$ and
\begin{equation}
\label{ineqn} \|F^{\lozenge n}\|_{p+2}^\prime\le
\frac{1}{A(2)}\left(A(2)\|F\|_p^\prime\right)^n,\quad
n=1,2,3\ldots
\end{equation}
\end{La}
{\bf Proof:} We proceed by induction. The case $n=1$ holds since
\[
\|F\|_{p+2}^\prime\le \|F\|_{p}^\prime,\quad {\rm for}\quad F\in
\mathcal H_p^\prime.
\]
Assume now that \eqref{ineqn} holds at rank $n$. Then, from
\eqref{vage} we have
\[
\begin{split}
\|F^{\lozenge(n+1)}\|_{p+2}^\prime&\le A(2)\|F\|_p^\prime
\|F^{\lozenge n} \|_{p+2}^\prime\\
&\le
A(2))\|F\|_p^\prime\frac{1}{A(2)}\left(A(2)\|F\|_p^\prime\right)^n\\
&=\frac{1}{A(2)}\left(A(2)\|F\|_p^\prime\right)^{n+1}.
\end{split}
\]
\section{The ring $\mathfrak{R}$}
\label{sec:r}
The Kondratiev space $S_{-1}$ endowed with the Wick product is a
commutative ring of sequences $(c_\alpha)_{\alpha\in\ell}$, with
properties $(P_1)$, $(P_2)$ and $(P_3)$, where in $(P_1)$ one
understands by evaluation at the origin the first coefficient of
the sequence. Using the Hermite transform (defined below), we
view $S_{-1}$ as a ring of powers series in infinitely many
variables. We point out that there are other commutative rings of
sequences with properties $(P_1)$, $(P_2)$ and $(P_3)$, and for
which a counterpart of the Hermite transform holds. See
\cite{a_salomon}.\\

The Hermite transform is defined by
\[
I(H_\alpha)=z^\alpha,\quad{\rm with}\,\, \alpha\in\ell\quad{\rm
and}\,\, z=(z_1,z_2,\,\ldots)\in\C^{\N}.
\]
Then
\[
I(H_\alpha\diamond H_\beta)=I(H_\alpha)I(H_\beta).
\]
It extends for $F=\sum_{\alpha\in\ell} a_\alpha H_\alpha\in
S_{-1}$ by the formula $I(F)(z)=\sum_{\alpha\in\ell} a_\alpha
z^\alpha$, and converges in sets of the form
\[
K_p(R)=\left\{z\in\mathbb C^{\mathbb N}\,\,: \,\, \sum_{\alpha\ne
0} |z|^\alpha (2\N)^{p\alpha}<R^2 \right\},
\]
where $p$ is such that $F\in\mathcal H_p^\prime$. The Kondratiev
space $S_{-1}$ is closed under the Wick product, and we have
\[
I(F\diamond G)(z)=I(F)(z)I(G)(z)\quad{\rm and}\quad
I(F+G)(z)=I(F)(z)+I(G)(z)
\]
for any $F,G\in S_{-1}$. Therefore the image of the Kondratiev
space $S_{-1}$ under the Hermite transform is a commutative ring,
denoted by
\[
\mathfrak{R}\Def\rm Im (I(S_{-1})).
\]
This ring was introduced in \cite{alp}. We transpose to it via the
Hermite transform the properties of $S_{-1}$. We have
\[
\mathfrak R=\bigcup_{k=1}^\infty I(\mathcal H_k^\prime).
\]

We define the adjoint $\mathbf{G}^*=(\mathbf{h}_{st})
\in\mathfrak R^{m\times n}$ of
$\mathbf{G}=(\mathbf{g}_{ts})\in\mathfrak R^{n\times m}$ by
$\mathbf{h}_{st}(z)=\mathbf{g}_{ts}^*(z)$ ($t\in\left\{1,\ldots,
n\right\}$ and $s\in\left\{1,\ldots, m\right\}$). Then for
$\mathbf{A}\in\mathfrak{R}^{n\times m}$ and $ \mathbf
{B}\in\mathfrak{R}^{m\times u}$ we have
\begin{equation}
(\mathbf{A}\mathbf{B})^*=\mathbf{B}^*\mathbf{A}^*.
\end{equation}
Note that $\mathbf{G}^*(0)=\mathbf{G}(0)^*$, where
$\mathbf{G}(0)^*$ is the usual adjoint matrix.

\begin{Dn}
An element $\mathbf{A}\in\mathfrak R^{n\times n}$ will be said
strictly positive, $\mathbf{A}>0$, if it can be written as
$\mathbf{A}=\mathbf{GG}^*$, where $\mathbf{G}\in \mathfrak
R^{n\times n}$ is invertible. It will be said positive if
$\mathbf{G}$ is not assumed to be invertible.
\end{Dn}

\begin{La}
Let $\mathbf{A}\in\mathfrak R^{n\times n}$. Then, $\mathbf{A}$ is
strictly positive if and only
if $\mathbf{A}(0)\in\mathbb C^{n\times n}$ is a strictly positive matrix
(in the usual sense).
\end{La}

{\bf Proof:} If $\mathbf{A}=\mathbf{GG}^*$ with $\det\mathbf{ G}(0)\not=0$,
then $\mathbf{A}(0)=\mathbf{G}(0)\mathbf{G}(0)^*$
is a strictly positive matrix. Conversely, assume that $\mathbf A\in
\mathfrak R^{n\times n}$ is such that $\mathbf A(0)>0$. We write
\[
\begin{split}
\mathbf A(z)&=\mathbf A(0)+(\mathbf A(z)-\mathbf A(0))\\
&=\sqrt{\mathbf A(0)}\{I_n+(
\sqrt{\mathbf A(0)})^{-1}(\mathbf A(z)-
\mathbf A(0))(\sqrt{\mathbf A(0)})^{-1}
\}\sqrt{\mathbf A(0)}.
\end{split}
\]
Let $\mathbf E(z)=\sqrt{\mathbf A(0)})^{-1}(\mathbf A(z)- \mathbf
A(0))(\sqrt{\mathbf A(0)})^{-1}$ vanishes at $z=(0,0,\ldots)$.
Property $(P_2)$ with
\[
f(\zeta)=1+\frac{1}{2}\zeta-\frac{1}{2\cdot 4}\zeta^2+\frac{1\cdot
3}{2\cdot 4\cdot 6} \zeta^3+\cdots=(1+\zeta)^{1/2},\quad |\zeta|<1,
\]
implies that
$
\mathbf A=\mathbf C^2,
$
where $\mathbf C=f(\mathbf E)=\mathbf C^*$.
\mbox{}\qed\mbox{}\\

Similarly, if $\mathbf A$ is positive, then $\mathbf A(0)$ is also positive, but the
converse statement need not hold. Take for instance $n=1$ and $\mathbf A(z)=z_1$. Then
it is readily seen that one cannot find $\rr\in\mathfrak R$ such that
$z_1=\rr^*(z)\rr(z)$.\\

We define the ring of polynomials with coefficients in
$\mathfrak{R}$ by $\mathfrak{R}[\lambda]$.
% where $\zeta\in\C$ (often we will allow in fact $\zeta\in\mathfrak{R}$).
To avoid confusion between the variable $\lambda$ and the
variables $z$ we introduce the notation
\[
\E(\rr)=\rr(0),\quad \rr\in\mathfrak R.
\]
\begin{Dn}
A rational function with values in $\mathfrak R^{n\times m}$ is an
expression of the form
\begin{equation}
\mathbf R(\lambda)=\mathbf p(\lambda) (\mathbf q(\lambda))^{-1}
\end{equation}
where $\mathbf p\in(\mathfrak{R}[\lambda])^{n\times m}$, and $
\mathbf q\in\mathfrak{R}[\lambda]$ and such that $\E(\mathbf
q(\lambda)) \not\equiv 0$.
\end{Dn}

Let $\mathbf R\in\mathfrak R^{n\times m}(\lambda)$. Then,
$\E(\mathbf R)\in\C^{n\times m}(\lambda)$, and it is readily seen
that
\begin{equation}
\label{eval} (\E(\mathbf R))(\lambda)=\E(\mathbf R(\lambda)).
\end{equation}

It is proved in \cite{alp} that every rational function with
values matrices with entries in $\mathfrak R$ and for which
$\E(\mathbf q(0))\not=0$ can be written as \eqref{real}.
\begin{Ex}
Let $\rr\in\mathfrak{R}$. The function
\[
F_\rr(\lambda)=(\lambda-\rr)(1-\lambda
\rr^*)^{-1}\in\mathfrak{R}(\lambda)
\]
is rational. It is defined for $\lambda\in\C$ such that
$1\not=\lambda (\E(\rr))^*$.
\end{Ex}
The next example of rational function need not be defined for
$\lambda=0$.
\begin{Ex}
Let $\rr\in\mathfrak{R}$. The function
\[
F_\rr(\lambda)=(\lambda-\rr)(\lambda-
\rr^*)^{-1}\in\mathfrak{R}(\lambda)
\]
is rational. It is defined for $\lambda\in\C$ such that
$\lambda\not=(\E(\rr))^*$.
\end{Ex}

\section{Analytic functions with values in $\mathfrak R$}
\setcounter{equation}{0} \label{sec:3} It is possible to define
analytic functions with values in a locally convex topological
vector space (see for instance the discussion in
\cite{grothendieck1, grothendieck2,
MR986066,herve}). Here the structure of
$\mathfrak R$ allows us to focus, locally, on the classical
definition of Hilbert space valued functions, as we now explain.

\begin{Pn}
Let  $\Omega\subset\C$ be an open set and let $\mathbf
f:\Omega\to\mathfrak{R}$  be a continuous function. Then,
$\mathbf f$ is locally Hilbert space valued, that is, for every
$\zeta_0\in\Omega$, there is a compact neighborhood $K$ of
$\zeta_0$ and a number $p_0$ such that $\mathbf f(K)\subset
I(\mathcal H^\prime_{p_0})$.
\end{Pn}

{\bf Proof:} Every $\zeta_0\in\Omega$ has a neighborhood $K$ of
the form $\overline{B_{\delta}}=\{\zeta\in\Omega~
;~|\zeta_0-\zeta|\leq\delta\}$ for some $\delta>0$. Since
$\bar{B_{\delta}}$ is a compact set and $\mathbf f$ is continuous,
$\mathbf f(\overline{B_{\delta}})$ is compact in $\mathfrak R$,
and therefore strongly bounded. See \cite[p. 54]{GS2_english}.
Thus there exists $p_0\in\N$ such that $\mathbf
f(\overline{B_{\delta}})\in I(\mathcal H^\prime_{p_0})$ and is
bounded in the norm of $I(\mathcal H^\prime_{p_0})$. See
\cite[Section 5.3 p. 45]{GS2_english}.
\mbox{}\qed\mbox{}\\

Therefore we can define an analytic function from $\Omega$ to
$\mathfrak R$ as a continuous function which locally admits a
power expansion with coefficients in one of the spaces
$I(\mathcal H^\prime_p)$. The following example shows that we
cannot expect to have a fixed $p$ in general.

\begin{Ex}
\label{ex:exex}
Let $\mathbf f(\lambda,z)=\sum_{n=1}^\infty
n^\frac{\lambda}{2}z_n$.
%H_{\epsilon^{(n)}}$.
Then $\mathbf f$ is continuous (as a function of $\lambda$) from
$\mathbb C$ into $\mathfrak R$, but there is no $p$ such that
$\mathbf f(\lambda,z)$ (viewed now as a function of $z$) belongs
to $I(\mathcal H^\prime_{p})$ for all $\lambda\in\C$.
\end{Ex}

Indeed, let $\lambda_0\in\mathbb C$. We have
\[
(\|\mathbf f(\lambda_0)\|^\prime_{p})^2=\sum_{n=1}^\infty
|n^{\lambda_0}| (2n)^{-p}=2^{-p}\sum_{n=1}^\infty n^{{\rm
Re}\,\lambda_0-p}<\infty,
\]
for $p>{\rm Re}\,\lambda_0+1$. To show continuity at a point
$\lambda_0\in\mathbb C$, we take $p>|\lambda_0|+2$, and restrict
$\lambda$ to be such that $|\lambda_0-\lambda|<1$. Using the
elementary estimate
\begin{equation}
\label{z1z2}
|e^{z_1}-e^{z_2}|\le |z_1-z_2|\cdot\max_{z\in[z_1,z_2]}|e^z|,
\end{equation}
for $z_1,z_2\in\mathbb C$, we have for $n=2,3,\ldots$
\[
|n^{\frac{\lambda}{2}}-n^{\frac{\lambda_0}{2}}|\le \ln n
\frac{|\lambda-\lambda_0|}{2} e^{\frac{|\lambda_0|+1}{2}\ln n}
\]
and so
\[
\begin{split}
(\|\mathbf f(\lambda)-\mathbf
f(\lambda_0)\|^\prime_{p+2})^2&=2^{-p-2}\sum_{n=2}^\infty
|n^{\frac{\lambda}{2}}-n^{\frac{\lambda_0}{2}}|^2n^{-p-2}\\
&\le\frac{|\lambda-\lambda_0|^2}{4}\sum_{n=2}^\infty\frac{(\ln
n)^2}{n^2} n^{|\lambda_0|+1-p}
\end{split}
\]
and hence the continuity at the point $\lambda_0$ in the norm
$\|\cdot\|^\prime_{p+2}$, and hence in $\mathfrak R$. See in particular \cite[p. 57]{GS2_english} for the latter.\\

Recall that, in the case of Hilbert space, weak and strong
analyticity are equivalent, and can be expressed in terms of
power series expansions. The argument uses the uniform
boundedness theorem. See \cite[Theorem VI.4, p. 189]{MR58:12429a}.
We define the evaluation of an $\mathfrak R$-valued analytic
function at a point $\rr\in\mathfrak{R}$. We first introduce
\[
{\mathfrak R}_\Omega=\{ \rr\in\mathfrak R;~\E(\rr)\in\Omega\},
\]
where $\Omega\subset\mathbb C$ is open.

\begin{Tm}
\label{tm:eval}
Let $\Omega$ be an open subset of $\mathbb C$,
and let $\mathbf f:\Omega\to\mathfrak{R}$ be an analytic
function. Let $\rr\in\mathfrak{R}_\Omega$, and let
\begin{equation}
\label{eqfzeta}
\mathbf f(\zeta)=\sum_{n=0}^\infty \mathbf
f_n(\zeta-\E(\rr))^n,
\end{equation}
be the Taylor expansion around $\E(\rr)\in\Omega$, where the
$\mathbf f_n\in\mathcal H_{p_0}^\prime$ for some $p_0\in\N$, and
where the convergence is in $\mathcal H_{p_0}^\prime$. The series
\begin{equation}
\label{frr}
\mathbf f(\rr)=\sum_{n=0}^\infty \mathbf
f_n(\rr-\E(\rr))^n
\end{equation}
converges in $\mathcal H^\prime_{q}$ for some $q>p_0$.
\end{Tm}

{\bf Proof:} Let $K$ be a compact neighborhood of $\E(\rr)$, and
let $p_0\in\mathbb N$ be such that $\mathbf f(K)\subset \mathcal
H_{p_0}^\prime$. Let furthermore $R$ be the radius of convergence
of the $\mathcal H_{p_0}^\prime$-valued power series
\eqref{eqfzeta}. In view of \eqref{lim=0}, there exists $p$,
which we can assume strictly larger than $p_0$, such that
\begin{equation}
\label{eq:radius}
A(2)\|\rr-\E(\rr)\|^\prime_{p}<R.
\end{equation}
On the other hand, using \eqref{ineqn}, we obtain
\[
\begin{split}
\|\mathbf f_n(\rr-\E(\rr))^{n}\|_{p+2}^\prime&\le A(2)\|\mathbf
f_n\|_{p_0}^\prime \|(\rr-\E(\rr))^{
n}\|_{p+2}^\prime\\
&\le \|\mathbf
f_n\|^\prime_{p_0}\left(A(2)\|\rr-\E(\rr)\|_{p}^\prime\right)^n.
\end{split}
\]
In view of \eqref{eq:radius}, the series
\[
\sum_{n=0}^\infty \|\mathbf
f_n\|^\prime_{p_0}\left(A(2)\|\rr-\E(\rr)\|_{p}^\prime\right)^n
\]
converges and so the series
\[
\sum_{n=0}^\infty \mathbf f_n(\rr-\E(\rr))^{n}
\]
converges absolutely in $I(\mathcal H_{p+2}^\prime)$.
\mbox{}\qed\mbox{}\\

The evaluation of $\mathbf f$ at $\rr$ is defined to be $\mathbf
f(\rr)$ given by \eqref{frr}.

\begin{Pn}
\label{pn:circ}
We can rewrite the evaluation at $\rr$ as a
Cauchy integral
\[
\mathbf f(\mathbf{\rr})=\frac{1}{2\pi i}\oint\frac{\mathbf
f(\zeta)}{\zeta- \rr}d\zeta
\]
where the integration is along a circle centered at $\E(\rr)$ and
of radius $r<R$ and in $\Omega$.
\end{Pn}
{\bf Proof:} As in Theorem \ref{tm:eval} we consider a compact
neighborhood $K$ of $\E(\rr)$, and let $p_0$ be such that that
$\mathbf f(K)\subset \mathcal H_{p_0}^\prime$. We consider a
simple closed path around $\E(\rr)$ which lies inside $K$.

We have
\[
\begin{split}
\frac{1}{2\pi i}\oint\frac{\mathbf f(\zeta)}{\zeta- \rr}d\zeta&=
\frac{1}{2\pi i}\oint\frac{\mathbf
f(\zeta)}{\zeta-\E(\rr)+\E(\rr)-
\rr}d\zeta\\
&=\frac{1}{2\pi i}\oint\frac{\mathbf f(\zeta)}{\zeta-
\E(\rr)}\left\{ \sum_{n=0}^\infty
\left(\frac{\mathbf{\rr}-\E(\mathbf{\rr})}{\zeta-\E(\mathbf{\rr})}\right)^n
\right\} d\zeta\\
&=\frac{1}{2\pi i}
\sum_{n=0}^\infty(\rr-\E(\rr))^n\Big\{\oint\frac{\mathbf
f(\zeta)}{(\zeta-\E(\rr))^{n+1}} d\zeta\Big\},
\end{split}
\]
where we have used the estimates as in the proof of Theorem
\ref{tm:eval} and the dominated convergence theorem to justify
the interchange of integration and summation.
\mbox{}\qed\mbox{}\\

Recall that a function $f$ analytic and contractive in the open
unit disk is called a Schur function. Furthermore, by the maximum
modulus principle, $f$ is in fact strictly contractive in $\mathbb
D$, unless it is identically equal to a unitary constant. We will
call a function $\mathbf f$ analytic from the open unit disk
$\mathbb D$ into $\mathfrak R$ a Schur function (notation:
$\mathbf f\in {S}_\mathfrak{R}$) if the function
\[
\lambda\mapsto\E(\mathbf f(\lambda))
\]
is a Schur function. For instance, both $1+z_1z_3$ and
$0.5+10z_1-3z_5$ are Schur functions. We now define the analog of
the open unit disk by
\[
\mathfrak{R}_\D= \{\rr\in\mathfrak{R};~\E(\rr)\in\D\},
\]
and the analog  of strictly contractive Schur functions as the
set of analytic functions $\mathbf f :\D\to\mathfrak{R}$ such
that the function $\lambda\mapsto \E(\mathbf f(\lambda))$ is a
strictly contractive Schur function.

\begin{Tm}
$\mathbf f\in{S}_\mathfrak{R}$ is a strictly contractive Schur
function if and only if $\mathbf f:\D\to\mathfrak{R}_\D$ is
analytic.
\end{Tm}
{\bf Proof:} If $\mathbf f$ is analytic from $\mathbb D$ into
$\mathfrak R$, and such that the
$\lambda\mapsto\E(\mathbf f(\lambda))$ is a strictly contractive
Schur function, it means by definition that the range of $\mathbf
f$ lies inside $\mathfrak{R}_\D$. Conversely, let $\mathbf
f:\D\to\mathfrak{R}$ be analytic and such that $\E(\mathbf f)$ is
a strictly contractive Schur function. Then for every $0<r<1$,
there exists $k\in\N$ (which may depend on $r$) such that $\mathbf
f(|\lambda|\leq r)\subset I(\mathcal{H^\prime}_{k})$. We can
write $\mathbf f$ as

\[
\mathbf f(\lambda)=\sum_{n=0}^\infty\lambda^n\mathbf{f_n},
\]
where $|\lambda|<r$ and $\mathbf{f_n}\in I(\mathcal H_k^\prime)$.
Now $\E(\mathbf
f)(\lambda)=\sum_{n=0}^\infty\lambda^n\E(\mathbf{f_n})$ for
$|\lambda|<r$. Since this holds for all $r\in(0,1)$ the function
$\lambda\mapsto\mathbf f(\lambda)$ has range inside
$\mathfrak{R}_\D$.\mbox{}\qed\mbox{}\\

\section{Nevanlinna-Pick Interpolation}
\label{sec:np}
In this section we solve the following interpolation problem
$(IP)$.
\begin{Pb}
\label{pb:paris_dec_2011} Given $n\in\N$ and points $\mathbf
a_1,\ldots, \mathbf a_n,\mathbf b_1,\ldots ,\mathbf
b_n\in{\mathfrak R}_\D$, find all Schur functions $\mathbf f$
with coefficients in $\mathfrak R$ such that $\mathbf f(\mathbf
a_i)=\mathbf b_i$ for $i=1,2,\ldots,n$.
\end{Pb}

The solution of this problem under the assumption that some
matrix is strictly positive, is presented in Theorem \ref{tm:np}
below. We first give some preliminary arguments, and note that if
$\mathbf f$ is a solution of the interpolation problem
\ref{pb:paris_dec_2011}, then $f=\E(\mathbf f)$ is a solution of
the classical interpolation problem
\begin{equation}
f(a_i)=b_i,\quad i=1,\ldots, n, \label{pbclassic}
\end{equation}
where we have set
\[
a_i=\E(\mathbf a_i)\quad{\rm and}\quad b_i=\E(\mathbf b_i),\quad
i=1,\ldots n.
\]
This last problem is solved as follows: let $P$ denote the
$n\times n$ Hermitian matrix with $ij$ entry equal to
\begin{equation}
\label{p} \frac{1-b_ib_j^*}{1-a_ia_j^*}.
\end{equation}
A necessary and sufficient condition for \eqref{pbclassic} to
have a solution in the family of Schur functions is that $P\ge
0$. We will assume $P>0$. Set, in the notation of the introduction,
\begin{equation}
\label{formulas_class}
\begin{split}
A_\zeta^*&=A={\rm diag}~(a_1^*,a_2^*\ldots ,
a_n^*),\\
-\begin{pmatrix}B_+\\ B_-\end{pmatrix}&=\begin{pmatrix}1&1&\cdots &1\\ b_1^*&b_2^*&\cdots
&b_n^*\end{pmatrix}\stackrel{\rm def}{=}C,\\
J&=\begin{pmatrix}1&0\\0&-1\end{pmatrix}.
\end{split}
\end{equation}
Furthermore, specializing the formula for $\Theta$ given in the introduction with $z_0=1$, or using the formula arising from
the theory of reproducing kernel Hilbert spaces (see  \cite{MR90g:47003},\cite{MR2002b:47144}), set
\[
\Theta(\lambda)=I_2-(1-\lambda)C(I_n-\lambda
A)^{-1}P^{-1}(I-A)^{-*}C^*J
\stackrel{\rm def.}{=}\begin{pmatrix}a(\lambda)&b(\lambda)\\
c(\lambda)& d(\lambda)\end{pmatrix}.
\]
We now gather the main properties of the matrix-valued function
$\Theta$ relevant to the present work. For proofs, we refer to
\cite{MR2002b:47144}, \cite{bgr}, \cite{MR90g:47003}.

\begin{Pn} The following hold:\\
$(a)$ The function $\Theta$ is $J$-inner with respect to the open unit
disk.\\
$(b)$ $\Theta$ has no poles in $\mathbb D$ and $a(\lambda)\sigma+b(\lambda)\not =0$
for all
$\lambda\in\mathbb D$ and all $\sigma$ in the closed unit disk.\\
$(c)$ The identity
\begin{equation}
\begin{pmatrix}1&-b_i\end{pmatrix}\Theta(a_i)=0,\quad i=1,\ldots
n. \label{paris_2011}
\end{equation}
is valid.
$(d)$ The linear fractional transformation
\[
T_{\Theta(\lambda)}(\sigma(\lambda))\stackrel{\rm
def.}{=}\frac{a(\lambda)\sigma(\lambda)+b(\lambda)}
{c(\lambda)\sigma(\lambda)+d(\lambda)}
\]
describes the set of all solutions of the problem
\eqref{pbclassic}  in the family of Schur functions when $\sigma$
varies in the family of Schur functions. \label{pn:theta}
\end{Pn}

  To solve the interpolation problem
\ref{pb:paris_dec_2011} we introduce the matrices $\mathbf
A,\mathbf C$ and $\mathbf P$, with entries in $\mathfrak R$,
built by formulas \eqref{p} and \eqref{formulas_class}, but with
$\mathbf a_1,\ldots, \mathbf a_n,\mathbf b_1,\ldots ,\mathbf b_n$
instead of $a_1,\ldots, a_n, b_1,\ldots , b_n$. Note that
$\mathbf P>0$ since $P>0$, and we can define the $\mathfrak
R^{2\times 2}$-valued function $\mathbf\Theta$ as $\Theta$ but
with $\mathbf A,\mathbf C$ and $\mathbf P$ instead of $A,C$ and
$P$. We have
\[
\E(\mathbf A)=A,\quad \E(\mathbf C)=C,\quad{\rm and}\quad
\E(\mathbf P)=P.
\]
Furthermore,
\begin{equation}
\E(\mathbf \Theta(\lambda))=\Theta(\lambda).
\end{equation}

\begin{Tm}
\label{tm:np}
Assume $\mathbf P>0$. Then, there is a one-to-one correspondence
between the solutions $\mathbf f$ of the problem
\ref{pb:paris_dec_2011} in ${S}_\mathfrak{R}$ and the elements
$\mathbf g\in{S}_\mathfrak{R}$ via the linear fractional
transformation $\mathbf f=T_{\mathbf \Theta}(\mathbf g)$.
\end{Tm}

{\bf Proof:} We first claim that the matrix-valued function
$\mathbf\Theta$ satisfies the counterparts of \eqref{paris_2011},
that is,
\begin{equation}
\begin{pmatrix}
1&-\mathbf b_i\end{pmatrix}\mathbf\Theta(\mathbf a_i)=0,\quad
i=1,\ldots n.
\label{paris_2011111}
\end{equation}
This is done using the permanence of algebraic identities. See
\cite[p. 456]{MR1129886} for the latter. Indeed, the  function
\[
(\det (I_n-\lambda A))(\det(I_n- A^*))(\det
P)(\prod_{\ell,j=1,\ldots n}(1- a_\ell a_j^*)) \Theta(\lambda)
\]
is a polynomial in $\lambda$ with coefficients which are
themselves polynomials in the $a_i$ and the $b_j$, with {\sl
entire} coefficients. Therefore, multiplying both sides of
\eqref{paris_2011} by the polynomial function
\[
(\det (I_n-\lambda A))(\det(I_n- A^*))(\det
P)(\prod_{\ell,j=1,\ldots n}(1- a_\ell a_j^*))
\]
evaluated at $\lambda=a_i$ ($i=1,2,\ldots, n$), and taking the
real and imaginary part of the equalities \eqref{paris_2011}, we
obtain for each $i$ four polynomial identities in the $4n$ real
variables ${\rm Re}\,a_j$, ${\rm Re}\, b_j$, ${\rm Im}\, a_j$,
${\rm Im}\, b_j$, with $j=1,\ldots n$, with {\sl entire}
coefficients, namely
\[
\begin{split}
{\rm Re}\,\left\{
\begin{pmatrix}1&-b_i\end{pmatrix}\det\,(I- a_iA)\det\,(I-A^*)\det\,P
\prod_{\ell,j=1,\ldots n}(1- a_\ell a_j^*)
\Theta(a_i)\right\}&=\\
&\hspace{-3cm}=\begin{pmatrix}0&0\end{pmatrix},\\
{\rm Im}\,\left\{
\begin{pmatrix}1&-b_i\end{pmatrix}\det\,(I- a_iA)\det\,(I-A^*)\det\,P
\prod_{\ell,j=1,\ldots n}(1- a_\ell a_j^*)
\Theta(a_i)\right\}&=\\
&\hspace{-3cm}=\begin{pmatrix}0&0\end{pmatrix}.
\end{split}
\]
It follows (see \cite[p. 456]{MR1129886}) that these identities
hold in any commutative rings, and in particular in $\mathfrak R$:
\[
\begin{split}
{\rm Re}\,\left\{
\begin{pmatrix}1&-\mathbf b_i\end{pmatrix}\det\,(I- \mathbf a_i
\mathbf A)\det\,(I-\mathbf A^*)\det\,\mathbf P
\prod_{\ell,j=1,\ldots n}(1- \mathbf a_\ell \mathbf a_j^*)
\mathbf \Theta(\mathbf a_i)\right\}&=\\
&\hspace{-3cm}=\begin{pmatrix}0&0\end{pmatrix},\\
{\rm Im}\,\left\{
\begin{pmatrix}1&-\mathbf b_i\end{pmatrix}\det\,(I- \mathbf a_i\mathbf A)\det\,(I-
\mathbf A^*)\det\,\mathbf P \prod_{\ell,j=1,\ldots n}(1- \mathbf
a_\ell \mathbf a_j^*)
\mathbf \Theta(\mathbf a_i)\right\}&=\\
&\hspace{-3cm}=\begin{pmatrix}0&0\end{pmatrix}.
\end{split}
\]
We now use the fact that we are in the ring $\mathfrak R$.
Because of the choice of the $\mathbf a_j$, the element
\[
\det\,(I- \mathbf a_i\mathbf A)\det\,(I- \mathbf
A^*)\prod_{\ell,j=1,\ldots n}(1- \mathbf a_\ell \mathbf a_j^*)
\]
is invertible in $\mathfrak R$. When furthermore $\mathbf P>0$ we
can divide both sides of the above equalities by
\[
\det\,(I- \mathbf a_i\mathbf A)\det\,(I- \mathbf
A^*)\prod_{\ell,j=1,\ldots n}(1- \mathbf a_\ell \mathbf
a_j^*)\det\,\mathbf P
\]
and obtain \eqref{paris_2011111}.\\

Let now $\rr\in S_{\mathfrak R}$, and let $\mathbf u,\mathbf v$
be analytic $\mathfrak R$-valued functions defined by
\[
\begin{pmatrix}\mathbf u(\lambda) \\ \mathbf v(\lambda)\end{pmatrix}=\mathbf \Theta
(\lambda)\begin{pmatrix}\rr (\lambda)\\
1\end{pmatrix}=\begin{pmatrix}\mathbf
a(\lambda)\rr(\lambda)+\mathbf b(\lambda)\\
\mathbf c(\lambda)\rr(\lambda)+\mathbf d(\lambda)\end{pmatrix}.
\]
Using \eqref{paris_2011111} we have that
\[
\mathbf u(\mathbf a_i)=\mathbf b_i\mathbf v(\mathbf a_i), \quad
i=1,\ldots, n.
\]
To conclude, we need to show that $\mathbf v(\mathbf a_i)$ is
invertible in $\mathfrak R$ for $i=1,\ldots n$. But we have
\[
\E(\mathbf v(\mathbf a_i))=c(a_i)\E(\rr)(a_i)+d(a_i),\quad
i=1,\ldots, n.
\]
Since the function $\Theta(\lambda)=\E(\mathbf \Theta(\lambda))$
is $J$-unitary on the unit circle and has no poles there.
Therefore, we have $c(a_i)\E(\rr)(a_i)+d(a_i)\not=0$ (see item
$(b)$ in Proposition \ref{pn:theta}), and hence $\mathbf v(\mathbf
a_i)$ is invertible in $\mathfrak R$. Therefore $\mathbf u\mathbf
v^{-1}=T_{\mathbf\Theta}(\rr)$ is a solution of
the interpolation problem.\\

Assume now that $\mathbf f$ is a solution. Then, we know from the
discussion before the theorem that there exists a Schur function
$\sigma(\lambda)$ such that
\begin{equation}
\label{rrr} \E(\mathbf f(\lambda))=T_{\E(\mathbf
\Theta(\lambda))}(\sigma(\lambda)).
\end{equation}
Define a $\mathfrak R$-valued function $\rr$ by
\[
\mathbf f(\lambda)=T_{\mathbf \Theta(\lambda)}(\rr(\lambda)).
\]
Taking $\E$ on both sides of this expression we obtain
\[
\E(\mathbf f(\lambda))=T_{\E(\mathbf
\Theta(\lambda))}(\E(\rr(\lambda))).
\]
Comparing with \eqref{rrr}, we obtain
$\E(\rr(\lambda))=\sigma(\lambda)$, and hence
$\rr\in{S}_\mathfrak{R}$.
\mbox{}\qed\mbox{}\\

\section{More general interpolation problem}
\setcounter{equation}{0}
The matrix-valued function $\Theta$ defined by \eqref{eq:theta} and describing the set of solutions of the bitangential
problem satisfies the conditions
\[
\begin{split}
\sum_{\lambda_0\in\mathbb D}{\rm
Res}_{\lambda=\lambda_0}
\begin{pmatrix}(\lambda I-A_\zeta)^{-1}B_+&B_-\end{pmatrix}\Theta(\lambda)&=0\\
\sum_{\lambda_0\in\mathbb D}{\rm
Res}_{\lambda=\lambda_0}\Theta(1/\lambda^*)^*\begin{pmatrix}C_-(\lambda I-A_\pi)^{-1}\\C_+\end{pmatrix}&=0.
\end{split}
\]
See also \cite{abds2}. As for the Nevanlinna-Pick case, these conditions can be translated into a finite number of
polynomial equations with coefficients in $\mathbb Z$, and the principle of permanence of identities
allows to extend these properties in the case of a commutative ring. On the other hand, we do not
know how to extend the third interpolation property, and so the
method is not  applicable to the most general bitangential interpolation
problem. On the other hand, if we restrict the parameter to be a
constant contractive matrix, the third condition also translates  into a
polynomial identity with entire coefficents, and the same method can
still be used.
The case of functions with poles inside the open unit disk, or the degenerate cases, are more difficult to
treat, and will be considered elsewhere.

\bibliographystyle{plain}
\def\cprime{$'$} \def\lfhook#1{\setbox0=\hbox{#1}{\ooalign{\hidewidth
  \lower1.5ex\hbox{'}\hidewidth\crcr\unhbox0}}} \def\cprime{$'$}
  \def\cprime{$'$} \def\cprime{$'$} \def\cprime{$'$} \def\cprime{$'$}

%\bibliography{/users/faculty/math/dany/bib/all}
%\bibliography{all}
%\bibliography{/home/alpay/bib/all}
%\bibliography{/home/dany/Desktop/bib/all}
\end{document}